\documentclass[AMS,Times1COL]{WileyNJDv5} 

\articletype{Article Type}%

\received{Date Month Year}
\revised{Date Month Year}
\accepted{Date Month Year}
\journal{Journal}
\volume{00}
\copyyear{2023}
\startpage{1}

\usepackage{xcolor}
\usepackage{caption}
\usepackage{subcaption}

\begin{document}

\title{Limiting behavior of mixed coherent systems with L\'evy-frailty Marshall-Olkin failure times}

\author[1]{Guido Lagos}
\author[1,2]{Javiera Barrera}
\author[3,4]{Pablo Romero}
\author[5]{Juan Valencia}

\authormark{LAGOS \textsc{et al.}}
\titlemark{LIMITING BEHAVIOR OF MIXED COHERENT SYSTEMS WITH L\'EVY-FRAILTY MARSHALL-OLKIN FAILURE TIMES}

\address[1]{\orgdiv{Facultad de Ingenier\'ia y Ciencias}, \orgname{Universidad Adolfo Ib\'a\~nez}, \orgaddress{\country{Chile}}}

\address[2]{\orgdiv{Department of Industrial Engineering and Management Sciences}, \orgname{Northwestern University}, \orgaddress{\city{Evanston}, \state{IL}, \country{USA}}}

\address[3]{\orgdiv{Facultad de Ingenier\'ia}, \orgname{Universidad de la Rep\'ublica}, \orgaddress{\city{Montevideo}, \country{Uruguay}}}

\address[4]{\orgdiv{Facultad de Ciencias Exactas y Naturales}, \orgname{Universidad de Buenos Aires}, \orgaddress{\city{Buenos Aires}, \country{Argentina}}}

\address[5]{\orgdiv{Departamento de Ingenier\'ia Industrial}, \orgname{Universidad de Santiago de Chile}, \orgaddress{\city{Santiago}, \country{Chile}}}

\corres{Corresponding author: Guido Lagos, Av Padre Hurtado 750 office A-225, Vi\~na del Mar, Chile. \email{guido.lagos@uai.cl}}



\abstract[Abstract]{In this paper we show a limit result for the reliability function of a system ---that is, the probability that the whole system is still operational after a certain given time--- when the number of components of the system grows to infinity.
More specifically, we consider a sequence of mixed coherent systems whose components are homogeneous and non-repairable, with failure-times governed by a L\'evy-frailty Marshall-Olkin (LFMO) distribution --- a distribution that allows simultaneous component failures.
We show that under integrability conditions the reliability function converges to the probability of a first-passage time of a L\'evy subordinator process.
To the best of our knowledge, this is the first result to tackle the asymptotic behavior of the reliability function as the number of components of the system grows.
To illustrate our approach, we give an example of a parametric family of reliability functions where the system failure time converges in distribution to an exponential random variable, and give computational experiments testing convergence.}


\maketitle

\newcommand{\I}[1]{ {\bf 1}_{\left\{ #1 \right\}} }
\newcommand{\RR}{\mathbb{R}}
\newcommand{\EE}{\mathbb{E}}
\newcommand{\PP}{ \mathbb{P}}
\newcommand{\dd}{\mathrm{d}}
\newcommand{\Qn}{Q_{n}}



\section{Introduction}

Computing and optimizing the reliability of a system is today an issue as relevant as ever, as our economy heavily relies on telecommunication networks, power grids and other complex networks.
Consequently, it is expected that they attain higher services levels~\cite{naldi2013towards,aceto2018comprehensive}, and to that end, societies are increasingly investing in extending and strengthening these networks.
Hence, it is essential to extend the study of reliability engineering ---a field in the literature since the 1960's--- by dropping technically convenient assumptions such as independent failures, as these can lead to overestimating the reliability; see~\cite{perez2018sixty} for further discussion.

In this paper we consider a system composed of a large number of components and study its \emph{reliability function}: a function giving the probability of the whole system being working by a certain given time.
This function depends on both the components' (random) lifetimes and their configuration, and from this perspective we consider the following setting.
Regarding component lifetimes, we assume that the components are non-repairable and fail at times distributed according to a \emph{L\'evy-frailty Marshall-Olkin} (LFMO) model~\cite{mai2009levy,barrera2020approximating}.
Regarding system configuration, we assume that the system is \emph{mixed coherent}~\cite[Ch.~2]{samaniego2007system}, where coherency means that more functional components cannot decrease the working probability of the system.

We consider the latter setting and give the first (to the best of the authors' knowledge) asymptotic result to make tractable the reliability function when the number of components of the system grows to infinity.
We give sufficient conditions under which the system failure time converges in distribution to the first-passage time of a L\'evy subordinator process.
Importantly, the latter are fundamental objects in the study of L\'evy processes and their fluctuations, see~\cite{kyprianou2014fluctuations}, and our result shows that first-passage times can also be seen as system failure times.
In particular, in the limit the system failure time is as amenable for estimation and simulation inasmuch as the L\'evy process is.

We consider the setting of \emph{coherent} systems, as coherency is an elementary property shared by many real-world systems~\cite[Ch.~2]{samaniego2007system}.
Indeed, it is an essential property in communication systems, transportation networks, logistics, risk management, decision sciences, among many others.
Nonetheless, there are exceptions where coherency is lost, such as in DC power grids~\cite{nesti2018emergent,guo2019less}  and some financial systems~\cite{capponi2016systemic}; the interested reader can consult~\cite{WangTrivedi} for the reliability analysis of non-coherent systems.
More specifically, we consider the setting of \emph{mixed} coherent systems~\cite[Ch.~3]{samaniego2007system}, which can be understood as the randomized choice among a finite collection of coherent systems (see also  \cite[Ch.~1]{navarro2021introduction}).
This additional source or randomness or ambiguity can aid in the modeling of complex systems with a large number of components, where the actual structure of the system may be misspecified or not precisely known.
In this setting the system can be parameterized by a (continuous) real vector, allowing to optimize the design using continuous optimization tools, whereas optimization over non-mixed coherent systems is usually intractable~\cite{lindqvist2021comparison}; see Section~\ref{subsec:mixed coherent} for the precise definition. 

Our work aims to contribute in the understanding of the capabilities and limitations of failure models that show dependence of the components' lifetimes.
In this sense, the Marshall-Olkin (MO) multivariate model, introduced in the 1960's~\cite{marshall1967multivariate}, is the cornerstone model for simultaneous failure of components; and still receives much attention in the communities of reliability and risk analysis~\cite{cherubini2015marshall}, as it is flexible enough to capture different settings, while having a rich mathematical structure to develop methodology.
Many attractive sub-families ---such as the LFMO model we consider here--- have been derived fairly recently, giving rise to lifetimes' models with a non-trivial dependency structure and limited complexity on the number of parameters, overcoming the exponential parametric complexity of a plain MO distribution~\cite{matthias2017simulating}.

Our results rely on the Samaniego signature, also known as the \emph{structural signature} of a system, that is related to the Barlow-Proschan index.
The former is a robust result that allows to efficiently compute the reliability of a general (coherent) system with exchangeable components~\cite{marichal2011signature}.
We pair this result with the LFMO model to estimate the reliability of a mixed coherent system with an increasing number of components.
Indeed, under the LFMO distribution the components are conditionally-iid and are governed by an underlying L\'evy subordinator that acts as a common latent factor  to all components.
In this way the parametric complexity of the model is bounded and limited to the parameters of the underlying L\'evy subordinator, and the addition of a new component only requires including a new standard exponential random variable to the model.

The hypotheses we make in our main result roughly require that, as the number of the components of the system grows, the probability that it continues to work even if a fraction of the components have failed has a well described limiting behavior.
More generally, our results can also be used to bound the reliability when there are (stochastic) bounds for the system signature.
Additionally, our experiments show that the approximation is valid for systems with only a moderate number of components. 

\subsection{Organization of this article}
A literature review is outlined in Section~\ref{sec:review}. The concepts of LFMO models and mixed coherent systems are presented in Section~\ref{sec:setting}.
To set the stage for our main result, in Section~\ref{sec:example} we give an illustrative application of the main theorem by giving a family of systems ruled by a single real parameter. 
Section~\ref{sec:main} states the main result of our paper,  Theorem~\ref{theo}, that gives a tractable asymptotic expression for the system reliability, and also gives a useful corollary.
Then, computational experiments testing convergence are given in Section~\ref{sec:experiments}.
The proof of the results shown are given in Section~\ref{sec:proofs}. 
Finally, Section~\ref{sec:conclusions} contains concluding remarks and trends for future work.

\subsection{Notation}

For $x$ in $\RR$ the value $\lceil x \rceil$ denotes the minimum integer $k$ such that $k \geq x$ and $\lfloor x \rfloor$ is the maximum integer $k$ such that $k \leq x$.
For an integer $n \geq 1$ the set $\Delta^n$ denotes the simplex of dimension $n$, i.e., $\Delta^n := \{p \in \RR^n : p \geq 0, \ \sum_i p_i = 1\}$.
For a vector $x$ in $\RR^n$, $\lvert x \rvert = \sum_{i=1}^n \lvert x_i \rvert$.
Lastly, $\#S$ denotes the cardinality of a set $S$.

\section{Literature Review}\label{sec:review}

Reliability Engineering is a long-standing field with a tradition of over 60 years of mathematical development.
It has received attention and essential contributions from several areas, such as Discrete Mathematics, Computer Science, Probability and Mechanical Engineering; see~\cite{2021-Brown}. Classical works in reliability deal with fixed-sized systems, whose closed-form expressions are usually involved; an overview of classical topics can be found in~\cite{hoyland2009system,gertsbakh2011network}.  To tackle essential quantities used in Reliability Engineering ---mean time-to-failure, mean replacement cost per cycle, among others--- it is usual to make simplified assumptions, such as iid lifetimes, parallel working components, and/or redundancy imposed by $k$-out-of-$n$  systems~\cite{eryilmaz2011estimation,nakagawa2012optimization}. However, as more real world systems rely on critical networks, they need to attain a higher service level; hence there is the need for more accurate reliability evaluation and optimization. Therefore, reliability modeling must improve to meet this challenge; for instance, consider models that add new features and move away from the classical assumptions that all components fail with iid probabilities.  The reader is invited to consult~\cite{perez2018sixty} for a recent historical survey and a discussion of these challenges. 

 Estimating the reliability of a network is a challenging problem. Even if arcs are assumed to fail as iid random variables, the computation of the source-terminal reliability and the all-to-all reliability are NP-complete problems~\cite{ball1986computational}. This problem has been undertaken by considering specific topologies such as a series-parallel, bipartite graphs or other underlying characteristics~\cite{satyanarayana1985linear,goharshady2020efficient}.  Others have proposed to tackle the problem using rare event simulation techniques (see \cite[Ch.~7]{rubino2009rare}). Many of these techniques allow us to consider independent failures that are not identically distributed without increasing their computational complexity.  The assumption that failures arise independently between components also needs to be left behind. It has been reported that components fail simultaneously more frequently than the independent assumption predicts~\cite{naldi2013towards}. Also, extreme weather events and natural disasters provoke geographically correlated failure, affecting the computation of reliability~\cite{rak2021disaster,vass2020model}.  In this context, the Marshall-Olkin (MO) distribution and other related models appear to be powerful tools for the modeling of components' lifetimes.
For instance, the MO model allows to capture geographical failure correlations~\cite{MatusIEEE}, it can be used to compute different reliability metrics~\cite{ozkut2019reliability}. Moreover, MO multivariate model, not only has exponentially distributed marginal life, it also has the memory lost property. This and other properties are used to simulate it efficiently  \cite{matthias2017simulating,Botev2015}. 

Asymptotic results for systems with many components can be helpful not only in giving closed-form expressions, but more importantly they can provide insights into the system. In this sense, the setting where all components have exchangeable failure times is well suited for asymptotic analysis, while also capturing the effect of dependent failures in the network reliability. To this end, we consider mixed coherent systems with non-repairable components and lifetimes given by the LFMO model. The latter was derived in~\cite{mai2009levy} as a parametric subfamily of an \emph{exchangeable} MO distribution, resulting in a \emph{conditionally-iid} distribution; see \cite[Ch.~3]{matthias2017simulating} for a comprehensive explanation of the MO subfamilies. In particular, it inherits the good properties of the MO distribution, a foundational tool in reliability modeling of simultaneous failures. The latter was originally proposed in the 1960's in~\cite{marshall1967multivariate}, and nowadays, there are entire international conferences dedicated to it; see, e.g.,~\cite{cherubini2015marshall}.

Another key ingredient of our main result is the use of the \emph{Samaniego signature} for mixed coherent systems~\cite[Ch.~3]{samaniego2007system}. This property was initially proposed in~\cite{samaniego1985closure} for iid failure times, and it has now been extended to exchangeable failure times~\cite{marichal2011signature, navarro2021conditions}, and even to indirect majority systems, such as some democratic voting systems~\cite{boland2001signatures} and non-exchangeable failure settings \cite{coolen2013generalizing}.
In short, the Samaniego result decomposes the reliability function by separating the system structure ---fully characterized in the system signature--- from the failure times of the components.
If two systems are equipped with the same components, the signature should be enough to reveal the most reliable one. Many influential studies were published during the last two decades boosting its study, see~\cite{naqvi2022system}. The computation of the system signature is a challenge by itself; see~\cite{yi2021signatures} for a comparison of methods for binary-state coherent systems. For example, the signature can be obtained for systems with a large number of components~\cite{da2018signature} by identifying subsystem's structures, even if there are shared components. These recent breakthroughs allow us to compute other important performance metrics such as the mean time to failure and the survival reliability using the signature~\cite{miziula2018bounds,eryilmaz2018mean,navarro2007reliability}. To estimate the survival reliability of the system with component lifetimes that are dependent random variables, the \emph{Samaniego signature} has been combined with copulas; see~\cite[Ch.~2]{navarro2021introduction} and  references therein.
A copula is a multivariate cumulative distribution function
of a random vector whose marginals are uniformly distributed on the interval $[0,1]$.
They are combined with the (generalized) inverse of any probability distribution to obtain a (dependent) multivariate vector by applying the Sklar's theorem; see \cite{matthias2017simulating}. In reliability, a natural assumption is that lifetimes are exponentially distributed, which is combined with classic copulas from the Archimedean family and with a signature result; see e.g.~\cite{eryilmaz2011estimation} and~\cite{torrado2021study}. As we already mentioned, MO models are natural candidates for lifetimes of dependent components, due to its characteristics. They have been paired with signatures to compute the reliability~\cite{bayramoglu2015reliability}, and to compute the mean time to failure in~\cite{bayramoglu2016mean}.
But even after assuming exchangeable components, the complexity of the expression allows to gain insight about the behavior of the system only through intensive computational efforts. In this context, the LFMO subfamily enables the asymptotic analysis of the reliability by, first, inheriting the good properties of the MO model and gaining other additional interesting properties, and second, by being a flexible model capable of capturing different dependence behavior between components' lifetimes.

Some exciting contributions based on the asymptotic analysis of reliability systems are~\cite{solovyev1991reliability}, where they estimate the probability of systemic failure when the repair rates, failure rates, and the number of components altogether behave asymptotically in a particular way.
Also,~\cite{bon1999approximation} extend the previous work to analyze certain settings for ageing properties of the components.
More recently,~\cite{bai2012shock} propose a cluster shock model for insurance risk with an underlying L\'evy process, and study the limit behavior when time goes to infinity; and~\cite{li2015network} perform a computational study for network reliability heuristically using arguments from Percolation Theory.
Perhaps the closest works to ours are~\cite{barrera2020limit}  and~\cite{navarro2008mean}. 
Indeed, the former article derives asymptotic limits for certain $k$-out-of-$n$  systems with underlying LFMO distribution, and the latter work derives an asymptotic expression for the reliability function and mean residual life of mixed coherent systems.
However, in both cases the asymptotic analysis is for large time and fixed-size systems, and in this paper we consider the opposite setting: fixed time and a growing number of components.
On a related vein, a fairly recent model that generalizes the LFMO distribution is~\cite{sun2017marshall}, where there is a multi-dimensional L\'evy subordinator that destroys  several components at once. This model was later applied for outages of energy power plants in \cite{malladi2020modeling}. Nonetheless, their approach is not as suitable for asymptotic analysis because the parametric complexity increases with the number of components of the  system. In fact, the multidimensional L\'evy subordinator has one dimension per component, which makes difficult to study the case when the number of components grows.

\section{Model}\label{sec:setting}

In this section we give the setting considered for our results.

\subsection{Mixed coherent systems}\label{subsec:mixed coherent}

In this section we address the property of \emph{mixed coherence} of a reliability system that can be working despite having components that are non-working --- say, a vehicle that can still work properly despite having some damaged fuses. 
The idea is to give a \emph{static} characterization of the system reliability, i.e., given a configuration of working and non-working components, 
to determine the probability that the whole system is still operational, irrespective of the time  instant when the non-working components failed. 
For that purpose, we focus on systems where more working components implies higher operational probability.

Consider a system with $n > 1$ components and assume the existence of a \emph{structure function} $\Phi : \{0,1\}^n \to \{0,1\}$, i.e., a deterministic function such that, for a given working/non-working configuration $x \in \{0,1\}^n$ of its components,  
where $x_i=1$ means that component $i$ is working and $0$ otherwise, then $\Phi(x) = 1$ means that the system is in working state and $0$ otherwise.

\begin{definition}[Coherent system]\label{def:coherent}
Consider a system with $n >1$ components and structure function $\Phi : \{0,1\}^n \to \{0,1\}$.
\begin{enumerate}
	\item Component $i$ is an \emph{irrelevant component} iff for all configurations $x \in \{0,1\}^n$, we have $\Phi( x \rvert_{x_i=0}) = \Phi(x \rvert_{x_i=1})$, where $ x \rvert_{x_i=y}$ denotes the configuration $x$ but with $x_i$ replaced by $y$.
	\item The system \emph{has no irrelevant components} iff none of its components is irrelevant.
	\item The system is \emph{monotone} iff for all component $i$ and all configuration $x$ in $\{0,1\}^n$ we have $\Phi( x \rvert_{x_i=0}) \leq \Phi( x \rvert_{x_i=1})$.
	\item The system is \emph{coherent} iff it is monotone and has no irrelevant components.
\end{enumerate}
\end{definition}

The following notion of \emph{signature} characterizes a coherent system.
It was introduced by Samaniego in~\cite{samaniego1985closure}, however here we give the definition of~\cite{marichal2011signature} based on~\cite{boland2001signatures}.

\begin{definition}[Signature]\label{def:signature}
Given a coherent system of size $n$ with structure function $\Phi$, we define its \emph{signature} $s$ in $\RR^n$ as
\begin{align*}
s_k := \phi_{n-k+1} - \phi_{n-k}, \qquad k = 1, \ldots, n,
\end{align*}
where 
\begin{align*}
\phi_{n-k} := & \sum_{x \in \{0,1\}^n, \lvert x \rvert=n-k} \Phi(x) \Biggm/ \binom{n}{n-k} 
\end{align*}
is the proportion of configurations, amongst the total number of configurations with exactly $k$ failed components, $\# \{x \in \{0,1\}^n : \lvert x \rvert=n-k \} = \binom{n}{n-k}$, that continue working despite having these $k$ failed components, $\# \{x \in \{0,1\}^n : \lvert x \rvert=n-k, \ \Phi(x)=1 \} = \sum_{x \in \{0,1\}^n, \lvert x \rvert=n-k} \Phi(x)$.
\end{definition}

Observe that because the system is monotone (since it is coherent) we have $0 = \phi_0 \leq \phi_1 \leq \ldots \leq \phi_n = 1$.
Hence, the signature $s$ of a system lies in the simplex of dimension $n$, i.e., in $\Delta^n := \{p \in \RR^n : p \geq 0, \ \sum_i p_i = 1\}$.

A noteworthy family \emph{coherent} systems is the $k$-out-of-$n$:F system (onwards denoted $k$-out-of-$n$ for simplicity), that fails at exactly the $k$-th failure of components.
In that case the signature of the system is the $k$-th canonical vector in $\RR^n$.
In particular, a series system has signature $(1, 0, \ldots, 0)$, and a parallel one $(0, \ldots, 0, 1)$.

We remark that the previous notions do not include any source of randomness; for instance, given a configuration of working and non-working components, the structure function gives a precise answer about the system being working or not.
However, we may want to randomize the latter, say because 
the complexity of the system makes it difficult to specify its precise structure, or especially if we want to optimize the system design over a continuous feasible set, as we will explain later.

The following notion of \emph{mixed} coherent system then arises as a way to include an additional source of randomness while trying to preserve some notion of coherence of the system.

\begin{definition}[Mixed coherent system and signature] \ 
\begin{enumerate}
	\item A system of size $n>1$ is \emph{mixed coherent} if it is the randomized choice, with certain given probabilities, amongst the set of coherent systems of size $n$.
That is, given the set of structure functions of coherent systems of size $n$, say $\{ \Phi^1 , \ldots, \Phi^m \}$, and given a discrete probability distribution $(p^1, \ldots, p^m)$ over the latter set, the structure function $\Phi$ of the mixed coherent system is random and equal to $\Phi^j$ with probability $p^j$, $j=1, \ldots, m$.
	\item For a mixed coherent system with structure function $\Phi$ its signature vector $s$ in $\RR^n$ is also $s_k := \phi_{n-k+1} - \phi_{n-k}$, for $k = 1, \ldots, n$, where $\phi_{n-k} := \sum_{x \in \{0,1\}^n, \lvert x \rvert=n-k} \EE \Phi(x) / \binom{n}{n-k}$ is now the \emph{expected} proportion of configurations that continue working despite having exactly $k$ failed components, amongst the total number of configurations that have exactly $k$ failed components.
\end{enumerate}
\end{definition}

Some comments are in order.
Note that any coherent system is also a mixed coherent one.
Also, the signature of a mixed coherent system is a vector in the simplex $\Delta^n$.
In fact, any choice of increasing sequence $0 = \phi_0 \leq \phi_1 \leq \ldots \leq \phi_n = 1$ defines the signature of a mixed coherent system of order $n$, whereas for non-mixed coherent ones there is only a finite number of such sequences (since there is a finite number of coherent systems).
A key observation that we exploit in our results is that any mixed coherent system, say with signature $s$, can be probabilistically interpreted as the randomized choice of the system $k$-out-of-$n$ with probability $s_k$, for each $k = 1, \ldots, n$.
Moreover, even if the system is coherent but not mixed, it still can be seen as a randomized choice between $k$-out-of-$n$ systems.
See~\cite[S.~2.3]{navarro2021introduction} for further details.

Perhaps the most appealing advantage of mixed coherent systems over non-mixed ones is the difficulty in imputing the signature of large and complex non-mixed system.
Indeed, it is known that the reliability estimation for a coherent system is computationally hard~\cite{1986Ball}.
This can be especially true for a practitioner modeling an involved real-world system.
In this case it may make more sense that the signature be specified by domain expert knowledge, by prescribing the values of the non-decreasing sequence $\phi$; as mentioned above, any such sequence does define a mixed coherent system, but may not define a non-mixed one.

Another attractive feature of mixed coherent systems ---over only coherent ones--- is that they allow to optimize the system design using continuous optimization tools, as they are parameterized by the system signature that lies in the simplex $\Delta^n$ of dimension $n$.
This is a salient advantage in the field of Reliability Economics, where the objective is to choose the system design that minimizes cost while meeting reliability constraints, or vice-versa (optimize risk under a budget constraint); see~\cite{2007SamaniegoNRL}.
Curiously, even in much simplified settings the optimal design can be attained in a strictly mixed coherent system; see, e.g.,~\cite[Example 7.1]{samaniego2007system}.
In contrast, optimization over (non-mixed) coherent systems can be a cumbersome task: the number of coherent systems grows exponentially with the size $n$, and finding the optimal coherent system is usually intractable~\cite[S.~7]{samaniego2007system}.
In fact, the sequence of number of monotone systems (see Definition~\ref{def:coherent}) of order $n$ is actually the Dedekind sequence, which is known only for $n < 10$, and asymptotic estimations are available~\cite{1991Wiedemann}. 

\subsection{L\'evy-frailty Marshall-Olkin distribution for failure times of components}\label{subsec:LFMO}

In this section we present the LFMO distribution.
It is a multidimensional distribution we use to model the failure times of the components of the system.

\begin{definition}\label{def:LFMO}
A random vector $T$ in $\RR^n$ is said to have a \emph{L\'evy-frailty Marshall-Olkin (LFMO)} distribution if its components $(T_1, \ldots, T_n)$ can be jointly defined as
\begin{align}\label{def eq:LFMO}
T_i := \inf \left\{ t \geq 0 \ : \ L_t > \varepsilon_i \right\}, \qquad i = 1, \ldots, n,
\end{align}
where $L = (L_t : t \geq 0)$ is a L\'evy subordinator stochastic process with $L_0=0$, and $\varepsilon_1, \ldots, \varepsilon_n$ is a collection of $n$ iid standard exponential random variables independent of $L$.
\end{definition}

\begin{figure}
{
\centering
\includegraphics[width=0.49\linewidth]{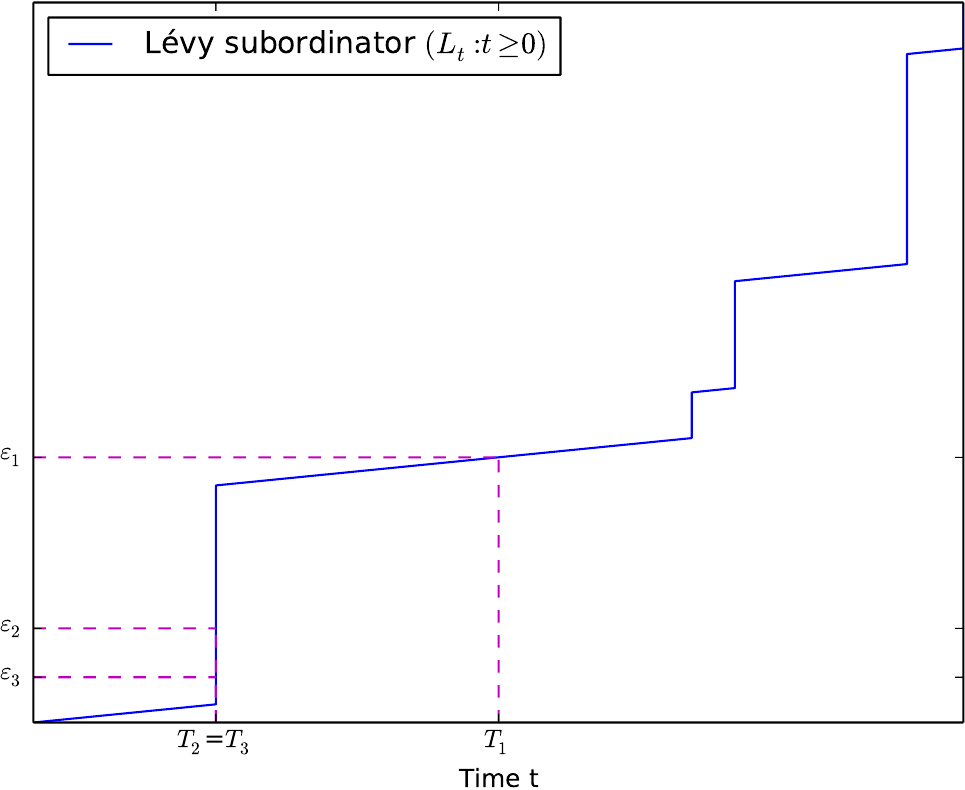}~
\includegraphics[width=0.49\linewidth]{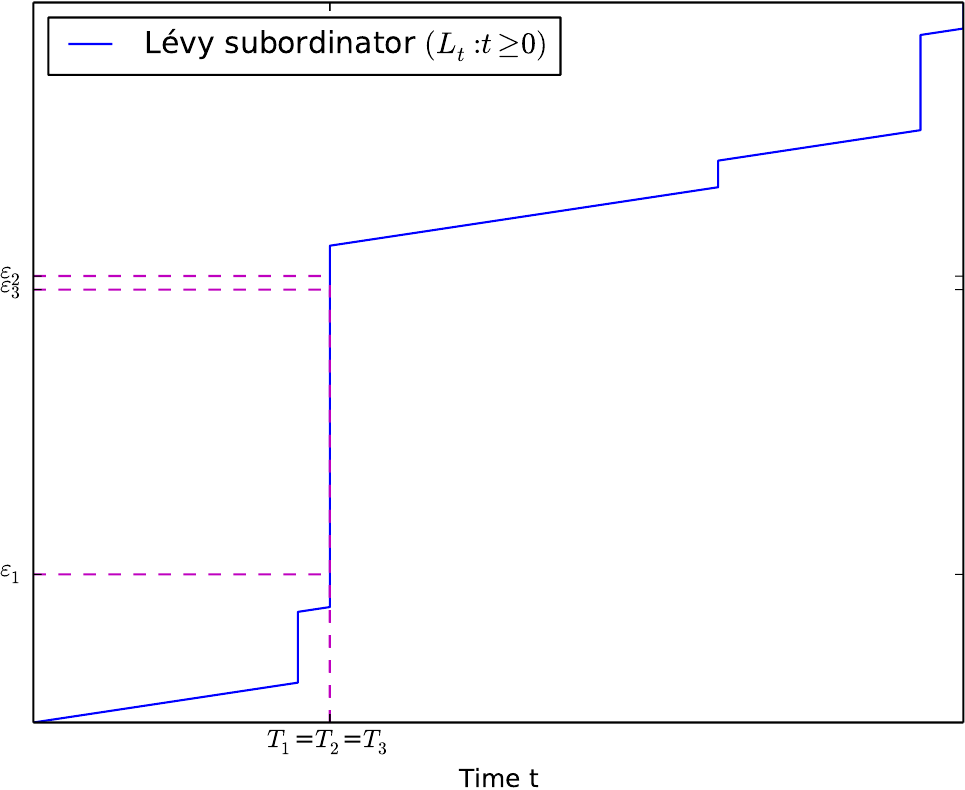}
\caption{\quad Two simulations of a random vector $(T_1, T_2, T_3)$ in $\mathbb{R}^3$ with an LFMO distribution: for each component~$i$, $T_i$ is the first time $t$ the L\'evy subordinator process $L$ surpasses the \emph{trigger} $\varepsilon_i$.
}
\label{fig:LFMO ex}
}
\end{figure}

The equation~\eqref{def eq:LFMO} gives a direct intuition of the model: component $i$ is ``destroyed'' when the subordinator $L$ up-crosses its so-called \emph{trigger} $\varepsilon_i$.
A possible modeling interpretation is that the components are homogenous in principle (hence the iid triggers), however are subject to a common stressor or deterioration force (the subordinator) that varies the speed at which the components fail and can make them fail in groups.
A pictorial example is the components of a vehicle subject to the stress of the mechanical vibrations.

We also remark that the LFMO distribution is especially useful for our modeling of a system growing in size: a trigger can be added at practically no parametric cost.

Importantly, the LFMO distribution can model from iid failures to simultaneous failure of components.
Indeed, on the one hand the deterministic pure-drift subordinator $L_t = \mu t$, for some $\mu>0$, produces failure times that are independent and exponentially distributed with rate $\mu$.
On the other hand, a jump in the subordinator can take down several component at once.
Figure~\ref{fig:LFMO ex} shows an example where the subordinator $S$ is a compound Poisson process with positive jumps and positive drift.
In general, though, the failure time for each individual component is marginally distributed as an exponential random variable with rate $\psi(1)$, where $\psi(x) := -\log \EE e^{-x L_1}$ is the Laplace exponent of the L\'evy subordinator $L$.
See~\cite{hering2012moment} for a moment-based estimation of the L\'evy subordinator associated to an LFMO distribution.

On a brief note about L\'evy subordinators, formally they are L\'evy processes with nondecreasing paths; see~\cite[Ch.~A.2]{matthias2017simulating} for a brief overview in the Marshall-Olkin (MO) distribution setting.
Intuitively speaking, though, they are continuous-time Markovian stochastic processes starting from zero, having non-decreasing sample paths, and can be understood as continuous-time versions of random walks with nonnegative jumps.
The simplest examples are the Poisson process and compound Poisson process, to which one can add a positive drift term (a non-decreasing deterministic linear term); and more involved examples are the running maximum of a Brownian motion, or the barrier-hitting time of a Brownian motion for each barrier value.
See~\cite{kyprianou2014fluctuations} for some nice lecture notes on L\'evy processes.

The LFMO distribution is a particular case of the MO distribution, proposed in~\cite{marshall1967multivariate}, and as such it inherits the interesting properties of the MO distribution: it generalizes to multiple dimensions the memoryless property of exponential random variables, and can be understood as a sequence of independent exponential ``shocks'' that takedown several components at once; see~\cite{barrera2020limit} and~\cite{bernhart2015survey} for further details.
The LFMO distribution was initially proposed in~\cite{mai2009levy} as a \emph{frailty} copula, and later in~\cite{mai2011reparameterizing} it was established as an \emph{extendible} version of a MO distribution with \emph{exchangeable} components, resulting in a \emph{conditionally-iid} distribution --- indeed, note that the components are iid when conditioned on the path of the L\'evy subordinator $L$.
See~\cite[Ch.~3]{matthias2017simulating} for a comprehensive explanation into these subfamilies and constructions.

\section{Asymptotically exponential system failure time}\label{sec:example}

To illustrate our main result, in this section we show Proposition~\ref{propos}, a direct corollary of  Theorem~\ref{theo}. The idea is that Proposition~\ref{propos} gives a glimpse of our result in a particularly interesting yet simple setting. Its proof is deferred to  Section~\ref{sec:proofs}.

\begin{proposition}\label{propos}
Consider a sequence of mixed coherent systems where the $n$-th system has $n$ components and signature vector $s^{(n)} = (s^{(n)}_1, \ldots, s^{(n)}_n)$ given by
\begin{align}\label{propos:sign}
s^{(n)}_k = C^{(n)} \left( n-k+1 \right)^{b-1} , \qquad k = 1, \ldots, n,
\end{align}
for some $b > 0$, and where $C^{(n)}$ is a constant that makes $\sum_{k=1}^n s^{(n)}_k = 1$.
Further, assume that the failure times of the system's components follow a L\'evy-frailty Marshall-Olkin distribution with underlying L\'evy subordinator $(L_t : t \geq 0 )$, common for all $n$.
If $T_\text{sys}^{(n)}$ the system failure time of the system with $n$ components, then we have
\begin{align}\label{lim:propos}
\PP \left( T_\text{sys}^{(n)} > t \right) \to  \PP \left( e_{\psi(b)} > t \right) \qquad \text{as } n \to  +\infty,
\end{align}
for all $t>0$, where $e_{\psi(b)}$ is an exponential random variable with rate $\psi(b)$ and $\psi(x) = -\log \EE e^{-x L_1}$ is the Laplace exponent of the L\'evy subordinator $L$.
In particular, the reliability function $R^{(n)}(t) = \PP ( T_\text{sys}^{(n)} > t )$ converges to $ e^{-\psi(b) t}$.
Moreover, in the limit the system fails when a random proportion $Q$ of the components have failed, where $Q$ has distribution $beta(1, b)$.
\end{proposition}

In simple words, Proposition~\ref{propos} states that, as the number of components increases, the distribution of the time at which the system fails approaches an exponential distribution with rate $\psi(b)$.
Interestingly, this gives a parametric family of systems where in the limit the failure time of the system is surprisingly simple.

An enginering question
is the qualification of the system signature given by Expression~\eqref{propos:sign}: when would one have such a system in hands?
Aside from systems specifically designed that way, there may also be cases where there is no accurate information on the system signature but where, nonetheless, there is a statistical notion on which failures tend to take down the system --- e.g., the first failure with some probability, the second failure with some other probability, and so on.
In that spirit, the parametric family of signatures~\eqref{propos:sign} can be interpreted as follows: the case $b=1$  corresponds to systems where each failure is equally likely to take down the system; the case $b$ in $(0,1)$ puts higher chance of failure when fewer components have failed (with higher probability of failure on fewer failed components as $b$ decreases to 0); and the case $b>1$ puts higher probability of  systemic failure as more components have failed (with higher chance on more failed components as $b$ increases to infinity).
This is reinforced by the fact that in the limit the system fails when e proportion $beta(1, b)$ of the components have failed, so in particular the mean proportion is $1/(1+b)$.
To illustrate, a series-parallel system with most components running in series will tend to fail when just a few components have failed, so possibly a $b$ close to 0 could be a better fit; whereas if most components run in parallel the system will tend to fail when a large number of its components have failed, so a larger $b$ could be a better fit.

Alternatively, for the signature $s^{(n)}$ with $s^{(n)}_k = C^{(n)} k^{b-1}$, $k = 1, \ldots, n$, for some $b>0$ and a proper normalizing constant $C^{(n)}$, it holds that for all $t>0$ the reliability function $R^{(n)}(t) = \PP ( T_\text{sys}^{(n)} > t )$ converges to $1 - \EE [ (1-e^{-L_t})^b ]$ as $n \to +\infty$.
The proof is analogous to the proof of Proposition~\ref{propos}, see Section~\ref{sec:proofs}.

\section{Main result}\label{sec:main}

In this section we give the main result of our work, Theorem~\ref{theo}, and give Corollary~\ref{corol}, a variation that expresses the hypotheses in terms of the system signature.

Consider a sequence of signatures $(s^{(n)} : n \geq 1)$, i.e., $s^{(n)}$ is in $\Delta^n = \{x \in \RR^n : x\geq 0, \ \sum_i x_i = 1\}$ for each $n\geq 1$.
Let $(Q_n: n \geq 1)$ be a sequence of discrete random variables where, for each $n$, $Q_n$ is supported on the set $\{\tfrac{1}{n+1}, \ \tfrac{2}{n+1}, \ \ldots, \ \tfrac{n}{n+1}\}$ and $\mathbb{P}(Q_n=k/(n+1))=s_{k}^{(n)}$ for each $k = 1, \ldots, n$.

\begin{description}

\item[\bf\hypertarget{hypA}{Hypothesis~(A)}] The sequence $( Q_n: n \geq 1 )$ converges in distribution, say to the random variable $Q$, i.e.,
$$Q_n \rightarrow_d Q.$$

\item[\bf\hypertarget{hypB}{Hypothesis~(B)}] The sequence $( Q_n: n \geq 1 )$ satisfies
\begin{align*}\label{theo bound}
\EE\left( \frac{1}{\sqrt{Q_n(1-Q_n)}} \right) = o(\sqrt{n}) \qquad \text{as } n \to  +\infty.
\end{align*}

\end{description}

We are now able to state the main result of our work.
Its proof is deferred to Section~\ref{sec:proofs}.

\begin{theorem}\label{theo}
Consider a sequence of mixed coherent systems whose components' failure times are distributed according to a L\'evy-frailty Marshall-Olkin distribution with common underlying L\'evy subordinator.
That is, for each $n \geq 1$ there is a mixed coherent system, say with signature vector $s^{(n)} = (s^{(n)}_1, \ldots, s^{(n)}_n)$ in $\Delta^n$, and the failure times of the $n$ components are jointly distributed according to a L\'evy-frailty Marshall-Olkin distribution over $\RR^n$, with underlying L\'evy-subordinator $(L_t : t \geq 0 )$ common for all $n$.
Denote by $T_\text{sys}^{(n)}$ the system failure time of the $n$-th system, and assume that the sequence of signatures $(s^{(n)} : n \geq 1)$ satisfies Hypotheses~\hyperlink{hypA}{(A)} and~\hyperlink{hypB}{(B)}.
Let also $Q$ be the limiting random variable of Hypothesis~\hyperlink{hypA}{(A)}.
Then $Q$ is independent of $L$, and for all $t>0$ we have
\begin{align}
\limsup_{n\to \infty} \Big\lvert \PP \left( T_\text{sys}^{(n)} > t \right) - \PP(\tau_{-\log(1-Q)} > t) \Big\rvert \leq \PP (L_t=-\log (1-Q))
\end{align}
where for a positive random variable $B$ we define $\tau_B := \inf \{ t > 0 : S_t > B \}$ as the first passage time of the subordinator $S$ across a value $B$.
In particular, if $L_t$ or $Q$ has a density then,
\begin{align*}
\PP \left( T_\text{sys}^{(n)} > t \right) \to
\PP(\tau_{-\log(1-Q)} > t) \qquad \text{as } n \to  +\infty
\end{align*}
i.e., $T_\text{sys}^{(n)}$ converges in distribution to $\tau_{-\log(1-Q)}$, and the reliability function $R^{(n)} (t) = \PP ( T_\text{sys}^{(n)} > t )$ converges to $\PP(L_t < -\log(1-Q))$ as $n \to  +\infty$.
\end{theorem}

We note that Hypotheses~\hyperlink{hypA}{(A)} and~\hyperlink{hypB}{(B)} make two requirements: (a)~that as the size of the system grows, the system fails essentially when a fraction of the total number of components have failed, and (b)~that the fractions where the system failure occurs stays sufficiently bounded away from 0 and from 1.
Indeed, for instance, the $k$-out-of-$n$ system (i.e., the system fails at exactly the $k$-th failure) for some \emph{fixed} $k$ and growing $n$ does satisfy (a) and~\hyperlink{hypA}{(A)}, but not~(b) nor~\hyperlink{hypB}{(B)}; whereas the system that fails when, e.g., 70\% of its components have failed does satisfy all requirements.
See~\cite{barrera2020approximating} for other more basic approximations for the following three cases of $k$-out-of-$n$ systems: that fail when just a few, a fraction, or most of the components have failed.

Another interesting case satisfying Theorem~\ref{theo} is the sequence of systems whose $n$-th system signature is given by $s^{(n)}_k = \binom{n-1}{k-1} p^{k-1} (1-p)^{n-k}$ for $k=1,\ldots,n$ and for a given $p$ in $(0,1)$, i.e., a binomial distribution of parameters $(n-1, p)$ shifted to be supported on $\{1, \ldots, n\}$.
In that case, by the weak law of large numbers, the limit distribution of Hypothesis~\hyperlink{hypA}{(A)} is the Dirac delta measure that puts all weight on the value $p$.
Hence, if $\PP(L_t = -\log(1-p)) = 0$, the reliability function $R^{(n)} (t) = \PP ( T_\text{sys}^{(n)} > t )$ converges to $\PP(L_t < -\log(1-p))$.

\subsection{Convergence of mean time-to-failure}

An important open question is when does the convergence of $T_\text{sys}^{(n)}$ to $\tau_{-\log(1-Q)}$ also holds in mean, i.e., the mean time-to-failure $\EE T_\text{sys}^{(n)} $ of the system converges to $\EE [ \tau_{-\log(1-Q)}] $.
As we show in the computational experiments in Section~\ref{sec:experiments}, there are cases where such convergence does not seem to hold, so there is the need to establish under which conditions it does.
This is an important question, as the mean time-to-failure $\EE T_\text{sys}^{(n)} $ is a crucial quantity in Reliability Engineering.

Furthermore, using~\cite[Proposition~1]{barrera2020limit} it can be shown that in the setting of Theorem~\ref{theo} the mean time-to-failure takes the following convoluted expression
\begin{align}
\EE \left[ T_\text{sys}^{(n)} \right] = \sum_{k=1}^n \sum_{l=n-k+1}^n \binom{n}{l} \binom{l-1}{n-k} \frac{s_k^{(n)}}{\psi(l)} (-1)^{l-n+k-1} .
\end{align}
where $\psi(x) = -\log \EE e^{-x L_1}$.
In contrast, to analyze and estimate the mean first-passage time $\EE [ \tau_{-\log(1-Q)}] $ of the subordinator $S$ we can use the literature on L\'evy processes and their simulation, see e.g.~\cite{asmussen2007stochastic,kyprianou2014fluctuations}.

\subsection{A corollary for explicit expressions for the signature}

We now give a corollary of Theorem~\ref{theo} for which the hypotheses may be easier to check if the expressions for the sequence of system signatures $s^{(n)}$ are readily available.

\begin{description}
	\item[\bf\hypertarget{hypC}{Hypothesis~(C)}] 
The sequence
 $(s^{(n)} : n \geq 1)$, with $s^{(n)}$ in $\Delta^n$ for each $n$, satisfies the pointwise convergence condition. 
\begin{align*}
n \, s^{(n)}_{\lceil n q \rceil} \to f(q) \quad \text{as $n \to +\infty$ \quad for all $q$ in $(0,1)$,}
\end{align*}
for some measurable function $f$ supported on the set $(0,1)$.
	\item[\bf\hypertarget{hypD}{Hypothesis~(D)}] 
	 $(s^{(n)} : n \geq 1)$
satisfies the domination condition
\begin{align*}
\sup_{n \geq 1} n \, s^{(n)}_{\lceil nq \rceil} \leq g(q) \quad \text{for all $q$ in $(0,1)$}
\end{align*}
for some measurable function $g$ satisfying $\int_0^1 g(q) \, \dd q < +\infty$.
\end{description}

\begin{corollary}\label{corol}
Consider the same setting of Theorem~\ref{theo} but where the sequence of signatures $(s^{(n)} : n \geq 1)$ instead satisfies Hypotheses~\hyperlink{hypC}{(C)} and~\hyperlink{hypD}{(D)}.
Then the function $f$ of Hypothesis~\hyperlink{hypC}{(C)} is the density of a random variable $Q$ supported on $(0,1)$, that is independent of $L$, and for all $t>0$ we have
\begin{align}
\PP \left( T_\text{sys}^{(n)} > t \right) \to
\PP(\tau_{-\log(1-Q)} > t) \qquad \text{as } n \to  +\infty.
\end{align}
In particular, $T_\text{sys}^{(n)}$ converges in distribution to $\tau_{-\log(1-Q)}$ and the reliability function $R^{(n)} (t) = \PP ( T_\text{sys}^{(n)} > t )$ converges to $\PP(L_t < -\log(1-Q))$ as $n \to  +\infty$.
\end{corollary}

The proof is deferred to Section~\ref{sec:proofs}, however it consists on showing that Hypotheses~\hyperlink{hypC}{(C)} and~\hyperlink{hypD}{(D)} together imply~\hyperlink{hypA}{(A)} and~\hyperlink{hypB}{(B)}.
The converse does not hold, as, for example, in the case of the system signatures
$s^{(n)}_k = \binom{n-1}{k-1} p^{k-1} (1-p)^{n-k}$ for $k=1,\ldots,n$ and fixed $p$ in $(0,1)$, we have that~\hyperlink{hypA}{(A)} and~\hyperlink{hypB}{(B)} hold, but neither~\hyperlink{hypC}{(C)} nor~\hyperlink{hypD}{(D)} hold.
Indeed, in this case the distribution of $Q_n$ converges to a Dirac delta measure supported on the value $p$; however $n \, s^{(n)}_{\lceil n q \rceil}$ tends to $0$ for $q \neq p$ and to $+\infty$ for $q=p$.

\section{Computational experiments}\label{sec:experiments}

In this section we show experiments that test computationally our results, and also the possible convergence in mean of the system failure time.

We test the convergence stated in Proposition~\ref{propos} for several values of the signature parameter $b$, and using as the L\'evy subordinator a compound Poisson process with nonnegative drift.
That is, our choice of $(L_t : t \geq 0)$ is
\begin{align}\label{def:CPP}
L_t = \mu t + \sum_{i=1}^{N_t} J_i
\end{align}
where $\mu \geq 0$ is a nonnegative \emph{drift} term, $(N_t : t \geq 0)$ is a Poisson process with rate $\lambda>0$, and $J_1, J_2, \ldots$ are nonnegative iid random variables that model the jumps of the subordinator, see e.g.~Figure~\ref{fig:LFMO ex}.
We choose such a process because any L\'evy subordinator can be approximated as close as desired by a compound Poisson processes, see~\cite[Ch.~XVII S.~2]{feller1971introduction}.
In this case we have $\EE L_t = (\mu + \lambda \EE J_1) t$ and the Laplace exponent $\psi$ of $L$ is $\psi(x) = \mu x + \lambda \left( 1 - \EE \left[ e^{-x J_1} \right] \right)$, see~\cite[Annex~A.2.1]{matthias2017simulating}.
We test the convergence of the system failure time $T_{sys}^{(n)}$ in Proposition~\ref{propos} for several parameters $\mu$, $\lambda$ and distributions for the jumps $J_i$ of the compound Poisson process, and for several values $b>0$ of the signature $s^{(n)}_k = C^{(n)} \left( n-k+1 \right)^{b-1}$ in Expression~\eqref{propos:sign}, where $C^{(n)} = 1/\sum_{i=1}^n i^{b-1}$ is a normalizing constant.

\subsection{Convergence in distribution}
To test the convergence in distribution~\eqref{lim:propos} of Proposition~\ref{propos} we would like to quantify how close the distribution of the system failure $T_{sys}^{(n)}$ is with respect to the limit exponential distribution of $e_{\psi(b)}$.
For that purpose, we use the p-value of the (two-sample) Kolmogorov-Smirnov test applied to two samples with 1,000 replications each of the random variables $T_{sys}^{(n)}$ and $e_{\psi(b)}$.
Here, a p-value close to 0, say less than 10\%, implies that statistically both samples come from different distributions.
Nonetheless, since the p-value of the test depends on the two samples and thus is essentially random, then we repeat 1,000 times each test and report the average p-value.
This procedure is used e.g.~in \cite{lachaud2006convergence} to test convergence of Markov chains, and in~\cite{barrera2020approximating} to test convergence in distribution of system failure times.

\begin{figure}[h]
{
\centering
     \begin{subfigure}[b]{.495\textwidth}
         \centering
         \includegraphics[width=\textwidth]{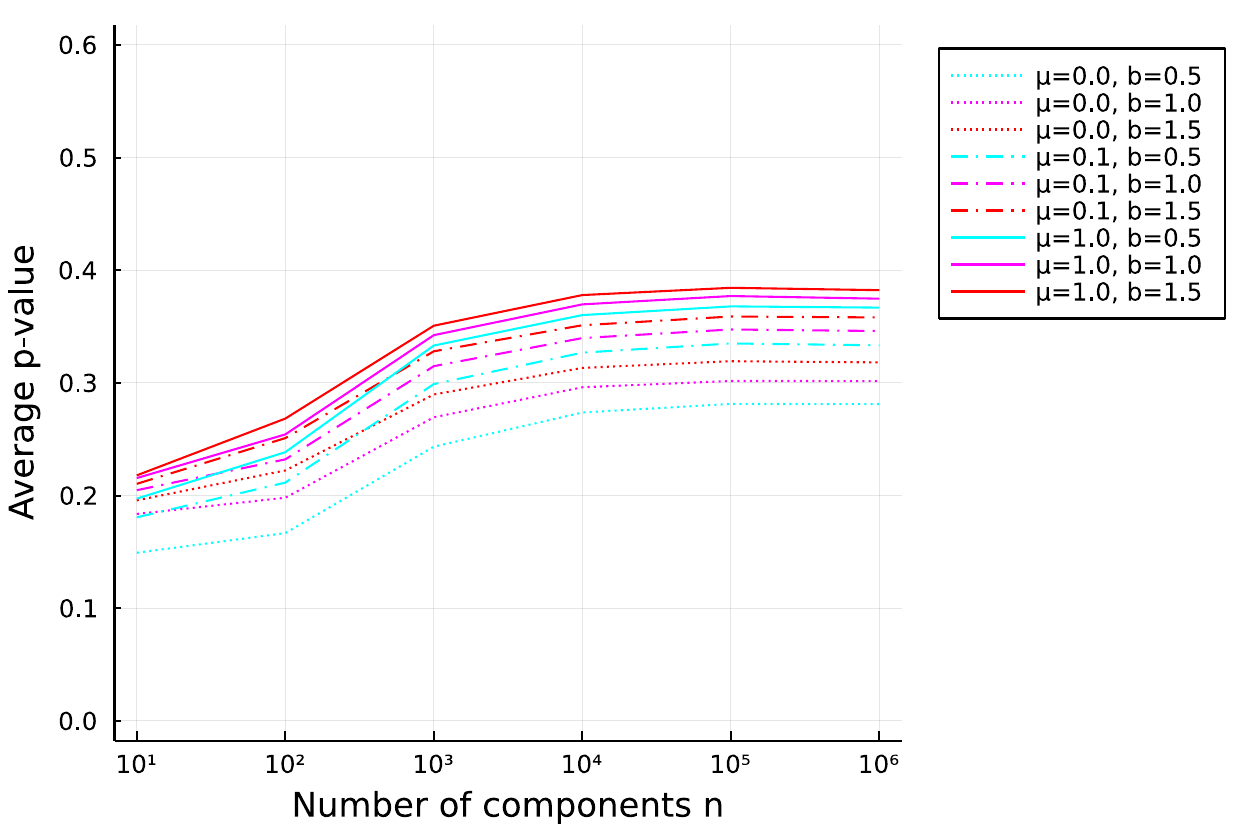}
         \caption{Uniform(0,1) jumps $J_i$}
         \label{fig:pvalues uniform}
     \end{subfigure}
     \hfill
     \begin{subfigure}[b]{.495\textwidth}
         \centering
         \includegraphics[width=\textwidth]{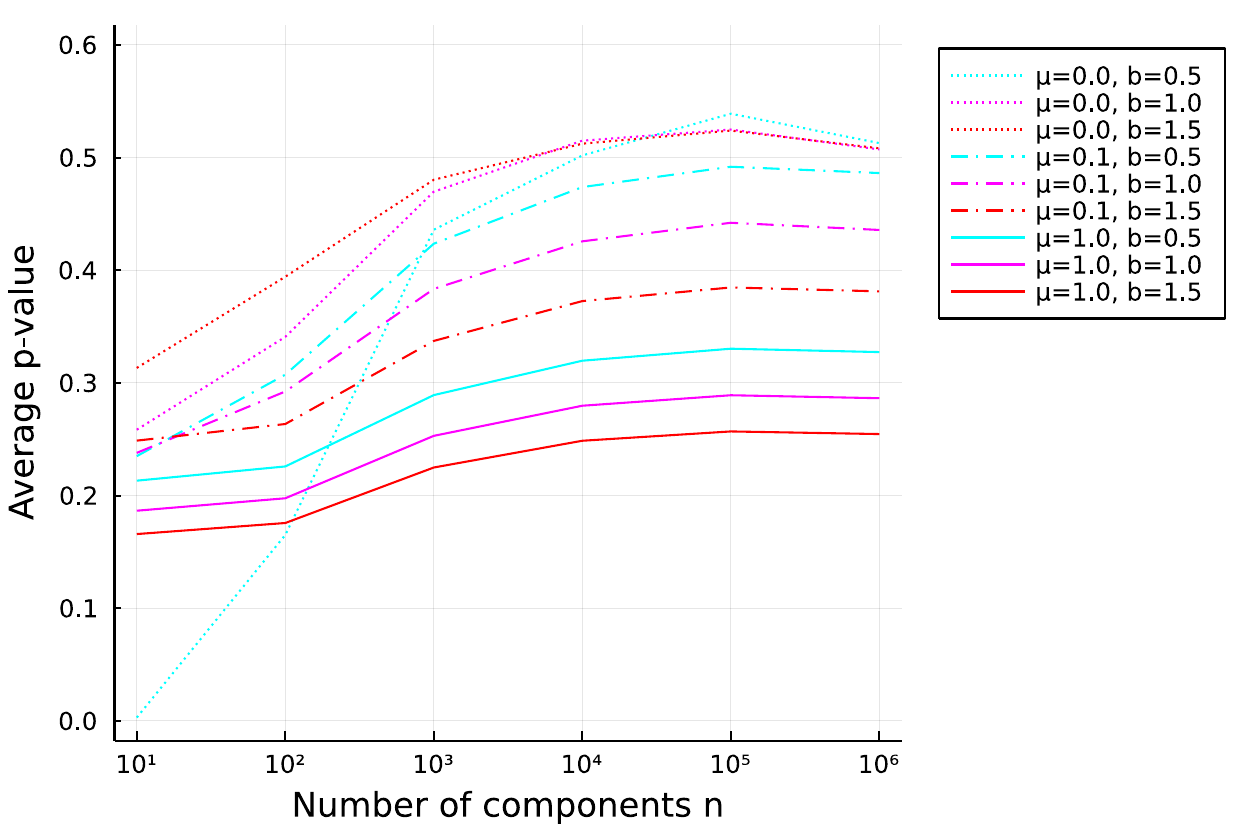}
         \caption{Exponential(1) jumps $J_i$}
         \label{fig:pvalues expo}
     \end{subfigure}
     \begin{subfigure}[b]{.495\textwidth}
         \centering
         \includegraphics[width=\textwidth]{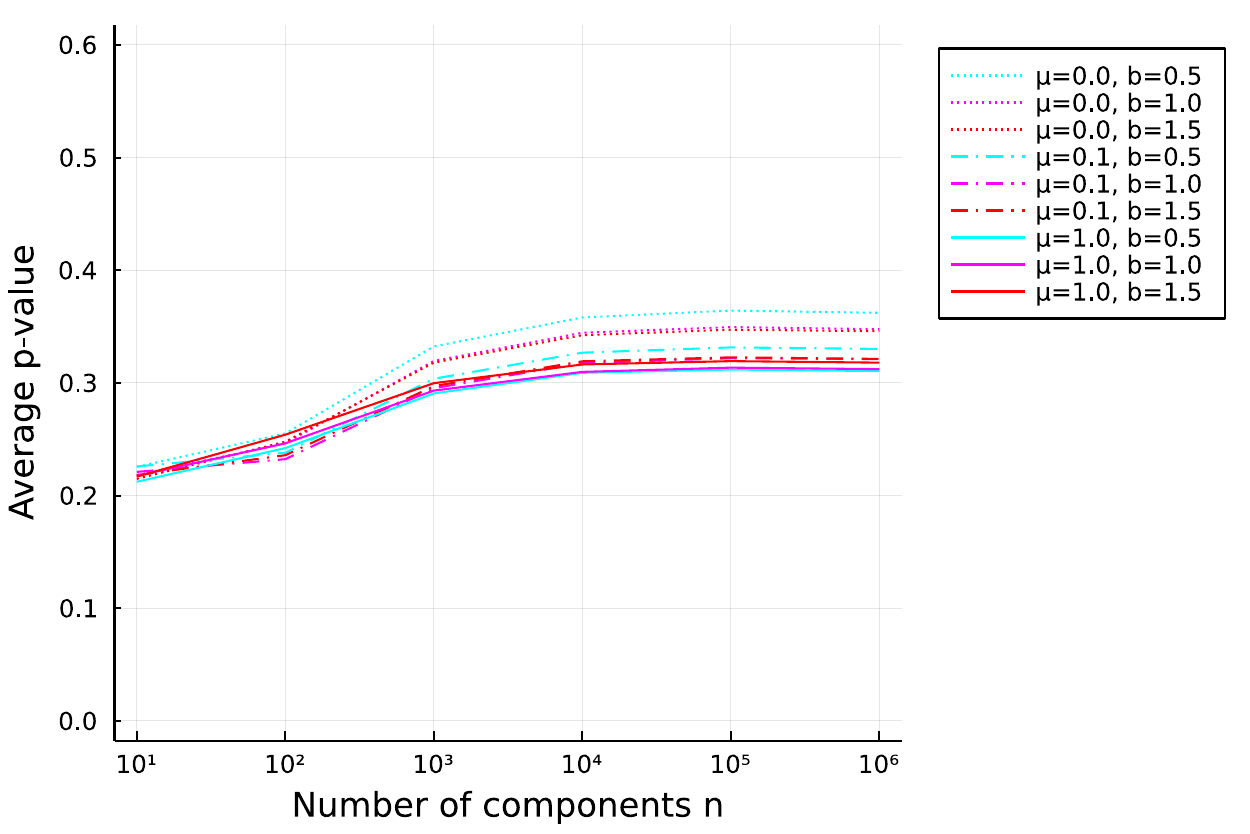}
         \caption{Pareto($\alpha=1.5$) jumps $J_i$}
         \label{fig:pvalues pareto}
     \end{subfigure}

\caption{
Average p-value over 1,000 repetitions of (two-sampled) Kolmogorv-Smirnov tests, comparing 1,000 samples of $T_\text{sys}^{(n)}$ and 1,0000 samples of the limit random variable $e_{\psi(b)}$, for several distributions of the jumps $J_i$ of the compound Poisson process subordinator~\eqref{def:CPP}.
For each case we show the values $\mu=0$, $0.5$ and $1$ for the drift term of the subordinator, and the values $b=0.5$, $1$ and $1.5$ for the parameter $b$ of the system signature in Expression~\eqref{propos:sign}.
}
\label{fig:pvalues}
}
\end{figure}

We consider several values for the parameters of the underlying L\'evy subordinator $L$ of the LFMO distribution, and for the parameters of the system signature.
Indeed, we consider several distributions for the jumps $J_i$ and the the drift value $\mu>0$ of compound Poisson process, and several values for the parameter $b>0$ of the system signature in~\eqref{propos:sign}.
We compute the aforementioned average p-value for several values of the number $n$ of components of the system, and focus on the case where $n$ grows exponentially: $n  = 10^1$, $10^2$, \ldots, $10^6$.
The results are shown in Figure~\ref{fig:pvalues}.

The experiments show that for all the test cases the limit in Proposition~\ref{propos:sign} provides a good approximation for as few as a hundred components.
Statistically speaking, for more than 100 components we cannot reject that both samples come from the same distribution.
Our experiments also suggest that the drift term can either facilitate or deter the convergence depending on the jump distribution.
Indeed, for uniform jumps, a positive drift term $\mu$  improves the convergence, whereas it is the opposite for exponential jumps.
An analogous behavior is observed with the parameter $b$: for uniform jumps, the convergence seems to be facilitated as $b$ grows, whereas it is the opposite for exponential jumps.
These findings suggest that a result on convergence speed should depend on these parameters.

\subsection{Convergence in mean}

\begin{figure}[h]
     \centering
     \begin{subfigure}[b]{0.495\textwidth}
         \centering
         \includegraphics[width=\textwidth]{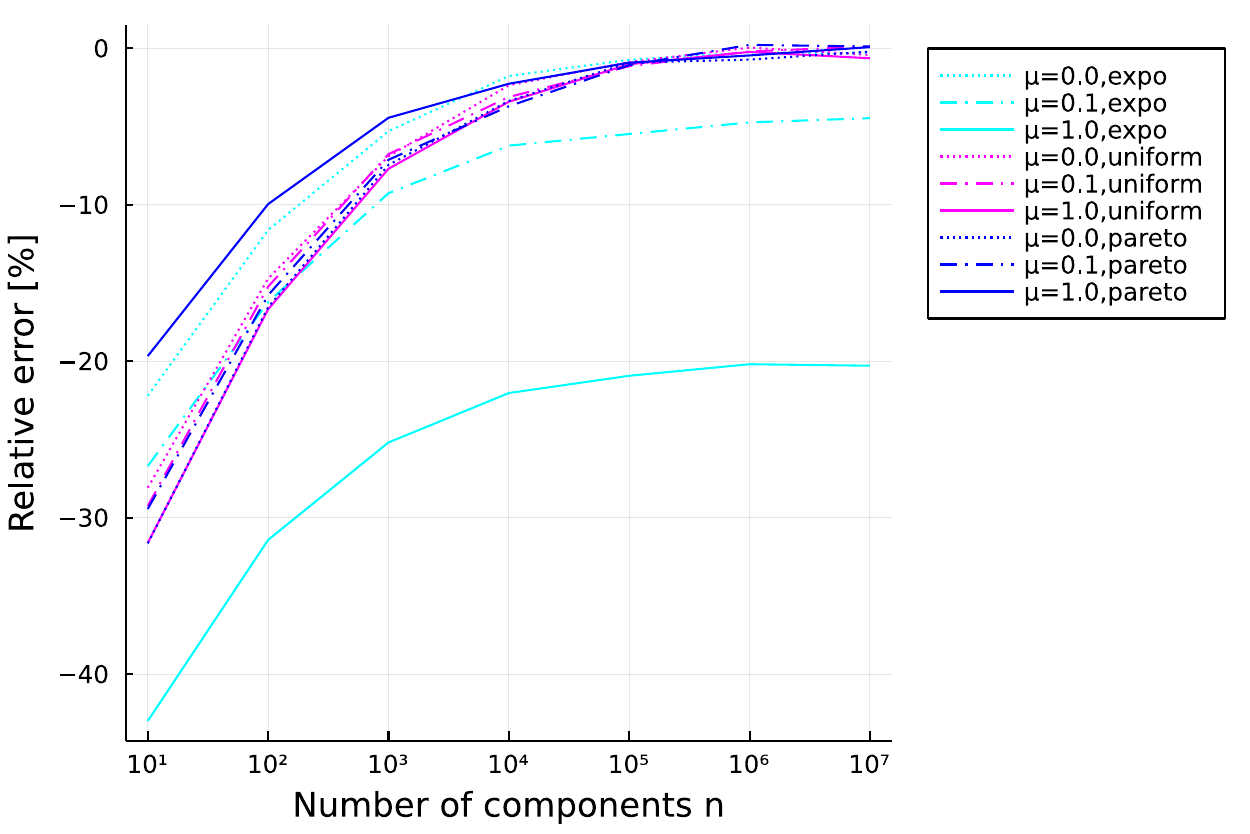}
         \caption{$b = 0.5$}
         \label{fig:means b05}
     \end{subfigure}
     \begin{subfigure}[b]{0.495\textwidth}
         \centering
         \includegraphics[width=\textwidth]{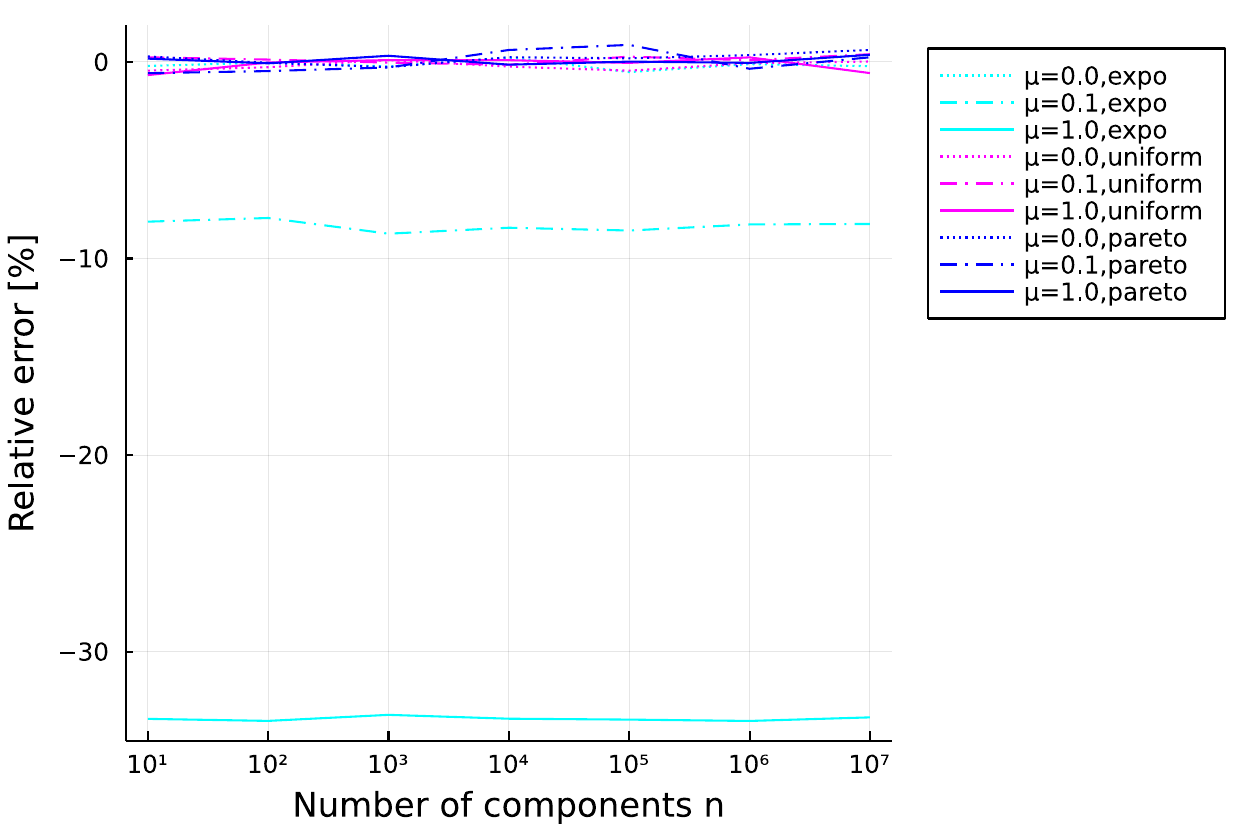}
         \caption{$b=1$}
         \label{fig:means b10}
     \end{subfigure}
     \begin{subfigure}[b]{0.495\textwidth}
         \centering
         \includegraphics[width=\textwidth]{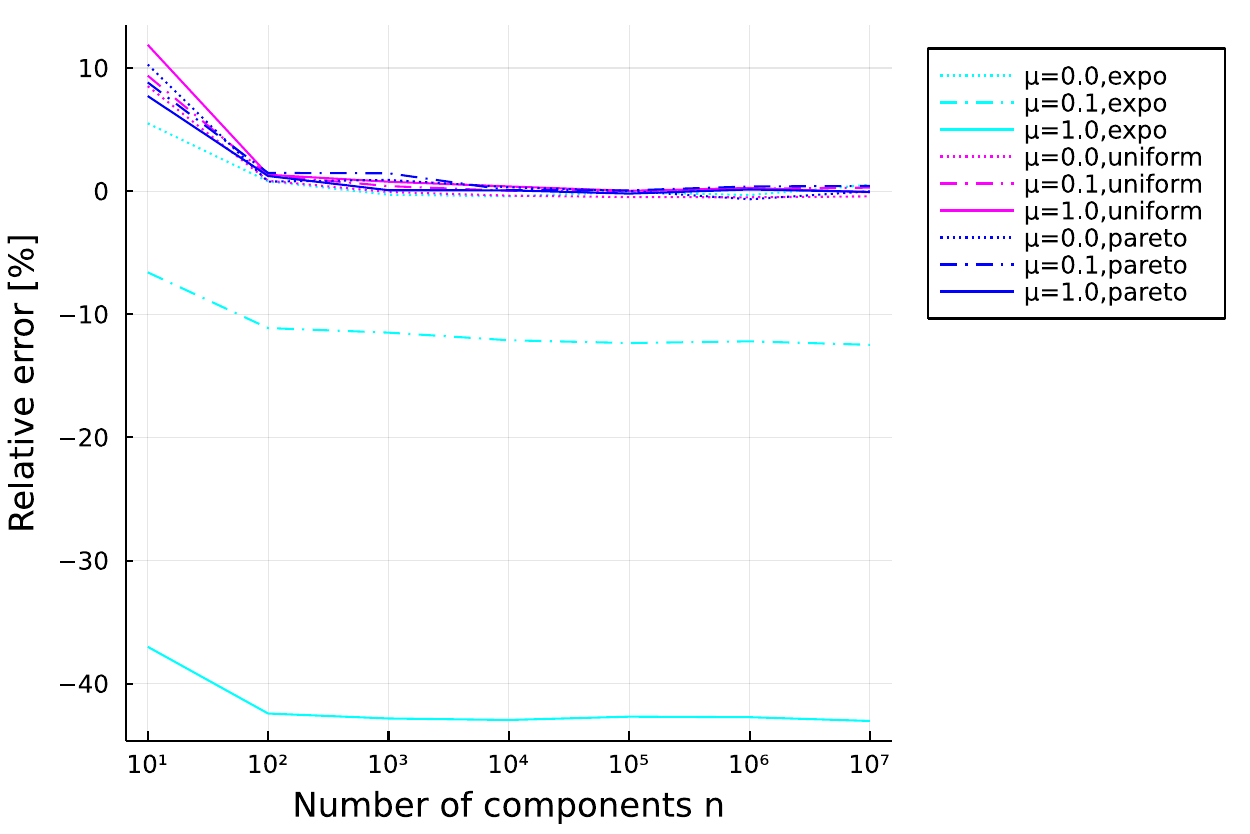}
         \caption{$b=1.5$}
         \label{fig:means b15}
     \end{subfigure}
\caption{
Relative error between the estimated mean (over 100,000 samples) of $T_\text{sys}^{(n)}$ and the mean $1/\psi(b)$ of the limit random variable $e_{\psi(b)}$ in~\eqref{lim:propos},
for several values of the parameter $b$ of the system signature in~\eqref{propos:sign}.
For each $b$ we consider the cases of the compound Poisson process subordinator~\eqref{def:CPP} having drift term $\mu=0$, $0.5$ and $1$, and jumps $J_i$ having an exponential(1), uniform(0, 1) and Pareto($\alpha=3$) distribution.
}
\label{fig:means}
\end{figure}

We now test numerically if the convergence~\eqref{lim:propos} stated in Proposition~\ref{propos} also holds in expectation, i.e., if, in the setting of the proposition, the mean $\EE T_\text{sys}^{(n)}$ converges to the mean $1/\psi(b)$ of the limit random variable $e_{\psi(b)}$.
For that, we take the average over 100,000 samples of $T_\text{sys}^{(n)}$, say $\overline{T}_\text{sys}^{(n)}$, and compute the relative error $( \overline{T}_\text{sys}^{(n)} - 1/\psi(b) ) \psi(b)$.
We perform these experiments for the L\'evy subordinator $L$ being the compound Poisson process~\eqref{def:CPP} with drift parameter $\mu=0$, $0.1$ and $1$; for jumps $J_i$ distributed as exponential(1), uniform(0,1), Pareto($\alpha=3$); and for the signature parameter $b$ taking values $b = 0.5$, $1$, and $1.5$.
Figure~\ref{fig:means} shows the results of these experiments.

We see that in all cases the mean seems to converge, except for exponential jumps with strictly positive drift term $\mu$.
Moreover, in this particular case the error between both means consistently grows as the drift increases.
In all other cases the error seems to converge to zero, with an initial negative bias when $b=0.5$, positive bias when $b=1.5$, and apparently no bias when $b=1$.

\section{Proofs}\label{sec:proofs}

In this section we give the proofs of Theorem~\ref{theo}, Corollary~\ref{corol} and Proposition~\ref{propos}.
The main ingredients to proving Theorem~\ref{theo} are applying the Samaniego signature result for distributions with exchangeable components, see~\cite{navarro2008mixture} and~\cite{samaniego2007system}, and then approximating the quantiles of iid exponential random variables in the spirit of~\cite[Proposition~2]{barrera2020approximating}.
The proof of Corollary~\ref{corol} consists on showing that Hypotheses~\hyperlink{hypC}{(C)} and~\hyperlink{hypD}{(D)} jointly imply Hypotheses~\hyperlink{hypA}{(A)} and~\hyperlink{hypB}{(B)}.
Lastly, Proposition~\ref{propos} comes from applying Corollary~\ref{corol} for the system signature given in~\eqref{propos:sign}.

\begin{proof}[Proof of Theorem~\ref{theo}.]
The proof starts by noting that for all $t > 0$ we have
\begin{align}
\nonumber
 \lefteqn{ \PP \left( T_\text{sys}^{(n)} > t \right)
= \sum_{k=1}^n \PP(T_{k:n} > t) \, s_k^{(n)}  } \\  \label{proof theo eq1}
&& = \sum_{k=1}^n \PP(T_{k:n} > t )\PP(\,(n+1)\Qn=k ) 
 =  \PP ( T_{ (n+1)Q_n : n} > t) .
\end{align}
Here, the first equality comes from the Samaniego signature for mixed coherent systems with exchangeable failure times, see~\cite{navarro2008mixture}, and since under the LFMO distribution the failure times $(T_i, i=1,\ldots,n)$ are exchangeable, see~\cite{mai2009levy} and~\cite[Ch.~3.2-.3]{matthias2017simulating}.
The second equality comes from $\PP(\Qn= k / (n+1) )= s_k^{(n)}$, by definition of $\Qn$.
We remark that $T_{ (n+1)\Qn : n}$ is always well defined, as $(n+1)\Qn \in \{1, \ldots, n\}$.

We will prove that
\begin{align*}
\limsup_{n \to \infty} \lvert \PP ( T_{(n+1) \Qn  : n} > t ) -  \PP (L_t<-\log (1-Q)) \rvert \leq \PP (L_t=-\log (1-Q))
\end{align*}
by bounding the difference of the left hand side as follows:
\begin{align*}
\lvert \PP ( T_{(n+1) \Qn  : n} > t ) -  \PP (L_t<-\log (1-Q)) \rvert \leq I_n + J_n,
\end{align*}
where $I_n$ and $J_n$ are defined as
\begin{align*}
I_n& =  \lvert \PP(T_{(n+1)\Qn : n} > t)- \PP\left(N_{0,1} > v(L_t,\Qn)\sqrt{n+1}\right) \rvert \\
J_n & =  \lvert  \PP\left(N_{0,1} > v(L_t,\Qn)\sqrt{n+1}\right)  -  \PP( L_t<-\log (1-Q) ) \rvert,
\end{align*}
where $v(l,q):= (l+\log(1-q))\sqrt{(1-q) / q}$, and  $N_{0,1}$ is a standard normal distribution independent of the subordinator $L$ and the sequence $(\Qn)_{n  \geq 1}$.
We will show that $I_n$ tends to zero and that $\limsup J_n \leq \PP(v(L_t,Q)=0) = \PP (L_t=-\log (1-Q))$.

We first prove that the term $I_n$ converges to zero.
For that purpose, the main tool will be a precise bound for the convergence of the order statistics of a collection of iid exponential random variables.
Recall that $(n+1)\Qn \in \{1, \ldots, n \}$ and note that by definition of the LFMO distribution we have
\begin{align}\label{proof theo eq2}
\PP(T_{(n+1)\Qn : n} > t) = \PP(L_t \leq \varepsilon_{(n+1)\Qn : n}),
\end{align}
where $\varepsilon_{1:n} \leq \ldots \leq \varepsilon_{n:n}$ are the increasing order statistics of the triggers $\varepsilon_1, \ldots, \varepsilon_n$ of the LFMO distribution.
We now use that the central order statistics of a sequence of iid exponential random variables converge in total variation to a normal distribution, with explicit bounds on the total variation.
Indeed, defining $q_n : = k/(n+1)$ for any fixed $k$ in $\{ 1, \ldots, n\}$, we have by~\cite[Corollary4.2.7]{reiss1989approximate} that the sequence $(\hat{\varepsilon}_{(n+1)q_n: n})_n$ defined as
\begin{align}\label{proof theo defqh} 
\hat{\varepsilon}_{(n+1)q_n: n} := \frac{\varepsilon_{(n+1)q_n: n} +\log(1-q_n)}{\sqrt{\frac{q_n}{(1-q_n)(n+1)}}},
\end{align}
converges in total variation to a standard normal distribution, and moreover
\begin{align}\label{proof theo lim2}
\| \PP_{\hat{\varepsilon}_{(n+1)q_n: n}}-N_{0,1} \|_\mathsf{TV} \leq \frac{C'}{\sqrt{q_n(1-q_n)(n+1)}}
\end{align}
for some constant $C'>0$, where $\| \cdot \|_\mathsf{TV}$ is the total variation distance.
In particular, using in~\eqref{proof theo lim2} the definition~\eqref{proof theo defqh} of $\hat{\varepsilon}_{(n+1)q_n: n}$, we obtain
\begin{align*}
& \sup_{l \in \RR} \Bigg\lvert \PP\left(\varepsilon_{(n+1)q_n : n}\geq l \right)-\PP\left( N_{0,1} > \frac{l+\log(1-q_n)}{\sqrt{\frac{q_n}{(1-q_n)(n+1)}}}\right) \Bigg\rvert \leq \frac{C'}{\sqrt{q_n(1-q_n)(n+1)}}.
\end{align*}
Therefore, using the definition of $I_n$ and~\eqref{proof theo eq2} we obtain that 
\begin{align*}
\lefteqn{ I_n \leq \EE\left[ \Big\lvert \PP\left(  \varepsilon_{(n+1)Q_n : n}\geq L_t \Big\vert L_t, \ Q_n \right)-\PP\left(  N_{0,1} > v(L_t,\Qn)\sqrt{n+1} \Big\vert L_t, \ Q_n \right) \Big\rvert \right] } \\
& = \EE\left[ \Bigg\lvert \PP\left( \varepsilon_{(n+1)Q_n : n}\geq L_t \Big\vert L_t, \ Q_n \right)-\PP\left( N_{0,1} > \frac{L_t+\log(1-Q_n)}{\sqrt{\frac{Q_n }{ (1-Q_n)(n+1)}}} \Bigg\vert L_t, \ Q_n \right) \Bigg\rvert \right] \\
&\qquad \leq \frac{1}{\sqrt{n+1}}\EE\left[ \frac{C'}{\sqrt{\Qn(1-\Qn)}} \right],
\end{align*}
which by Hypothesis~\hyperlink{hypB}{(B)} vanishes as $n$ goes to $+\infty$.
This shows that $I_n$ tends to $0$.

We will now show that $\limsup J_n \leq \PP(L_t=- \log(1-Q))$.
Note first that
\begin{align}
J_n &= \lvert \PP(N_{0,1} \leq v(L_t,Q_n)\sqrt{n+1} )-\PP(L_t \geq -\log(1-Q)) \rvert \nonumber \\
 &= \lvert \EE(\varphi( v(L_t,Q_n)\sqrt{n+1}) )-\PP(v(L_t,Q) \geq 0) \rvert, \label{eq:theo proof}
 \end{align}
where, recall, $v(l, q) = (l+\log(1-q))\sqrt{(1-q) / q}$, and $\varphi$ is the cumulative distribution function of a standard normal random variable.

Now, by Hypothesis~\hyperlink{hypA}{(A)}, $Q_n$ converges in distribution to $Q$, so in particular the random tuple $(L_t, Q_n)$ also converges in distribution to $(L_t, Q)$.
Skorokhod's representation theorem~\cite[S.~6]{billingsley1999convergence} guarantees the existence of a probability space, say with probability measure $\widetilde\PP$, where the sequence $( (L_t, Q_n) : n\geq 1 )$ converges $\widetilde\PP$-a.s.~to $(L_t, Q)$; this holds if the topological space where the limit random variable takes values is separable, which is the case for $\RR_+ \times (0, 1)$.
Furthermore, as $v$ is continuous on the latter domain then $v(L_t, Q_n)$ converges $\widetilde\PP$-a.s.~to $v(L_t,Q)$.
Then, since $\varphi$ is continuous, for $v(L_t,Q) \neq 0$ we have $\widetilde\PP$-a.s.~that $\lim_n \varphi(v(L_t,Q_n)\sqrt{n+1}) = 1$ if $v(L_t,Q)>0$ and $=0$ if $v(L_t,Q)<0$.
So in particular, $\widetilde\PP$-a.s.,
\begin{align*} 
\liminf_{n \to +\infty } \varphi(v(L_t,Q_n)\sqrt{n+1}) & \geq
\I{v(L_t,Q)>0} \\
\limsup_{n \to +\infty } \varphi(v(L_t,Q_n)\sqrt{n+1})& \leq   \I{v(L_t,Q)\geq 0}.
\end{align*}
Hence, using Fatou's lemma and its reverse we obtain
\begin{align}
\liminf_{n \to +\infty}  \EE(\varphi( v(L_t,Q_n)\sqrt{n+1} )) \label{proof ineq 1}
& \geq \PP(v(L_t,Q) > 0), \\
\limsup_{n \to +\infty}  \EE(\varphi( v(L_t,Q_n)\sqrt{n+1} )) & \leq \PP(v(L_t,Q)\geq 0). \label{proof ineq 2}
\end{align}

Finally, plugging inequalities~\eqref{proof ineq 1} and~\eqref{proof ineq 2} into~\eqref{eq:theo proof} we obtain
 \begin{align*}
\limsup J_n
\leq  \PP(v(L_t,Q)= 0) =  \PP(L_t=- \log(1-Q)).
\end{align*}
Furthermore, if $Q$ or $L_t$ have a density then $\PP(L_t=- \log(1-Q))=0$.
This concludes the proof of Theorem~\ref{theo}.
\end{proof}

We now give the proof of Corollary~\ref{corol}.

\begin{proof}[Proof of Corollary~\ref{corol}.]
Assuming that Hypotheses~\hyperlink{hypC}{(C)} and~\hyperlink{hypD}{(D)} hold, we will prove that Hypotheses~\hyperlink{hypA}{(A)} and~\hyperlink{hypB}{(B)} are satisfied. The result then follows from Theorem~\ref{theo}.

To prove~\hyperlink{hypA}{(A)} define first $f_n(q) := n \, s^{(n)}_{\lceil nq \rceil} $ for $q \in (0,1)$ and note that
\begin{align*}
\PP(Q_n	\leq q) = \sum_{k=1}^{\lfloor (n+1) q \rfloor} s^{(n)}_k = \sum_{k=1}^{\lfloor (n+1) q \rfloor} \int_{\tfrac{k-1}{n}}^{\tfrac{k}{n}} f_n(x) \dd x = \int_0^1 f_n(x)\I{n x \leq  \lfloor (n+1) q \rfloor } \dd x.
\end{align*}
Note too that $\lfloor (n+1) q \rfloor / n \to q$ as $n \to +\infty$, since $ q-\frac{1-q}{n}\leq \frac{ \lfloor (n+1) q \rfloor }{n}  \leq  q+\frac{q}{n} $ for all $q$ in $(0,1)$.
Then, using the dominated convergence theorem with Hypotheses~\hyperlink{hypC}{(C)} and~\hyperlink{hypD}{(D)} we obtain that $\lim_{n \to \infty} \PP(Q_n \leq q) = \int_0^q f(x)~\dd x$ for the measurable function $f$ such that $\int_{0}^{1}f(x)dx < +\infty$, i.e.,~Hypothesis~\hyperlink{hypA}{(A)} holds. 

To prove~\hyperlink{hypB}{(B)} we first decompose
\begin{align*}
\EE\left[ \frac{1}{\sqrt{Q_n(1-Q_n)}} \right] & = \EE\left[ \frac{\I{Q_n \in [\frac{1}{\sqrt{n}}, 1-\frac{1}{\sqrt{n}}]}}{\sqrt{Q_n(1-Q_n)}} \right] + \EE\left[ \frac{\I{Q_n \notin [\frac{1}{\sqrt{n}}, 1-\frac{1}{\sqrt{n}}]}}{\sqrt{Q_n(1-Q_n)}} \right].
\end{align*}
Now, note that for $\delta$ in $(0, 1/2)$ we have
\begin{align*}
\frac{\I{q \in [\delta, 1-\delta]} }{ \sqrt{q (1-q)} } \leq \frac{1}{ \sqrt{\delta(1-\delta)} } \leq \frac{1}{\sqrt{\delta}}.
\end{align*}
The latter implies that for $n \geq 5$ it holds
\begin{align*}
\EE\left[ \frac{\I{Q_n \in [\frac{1}{\sqrt{n}}, 1-\frac{1}{\sqrt{n}}]}}{\sqrt{Q_n(1-Q_n)}} \right] \leq \EE\left[ \frac{\I{Q_n \in [\frac{1}{\sqrt{n}}, 1-\frac{1}{\sqrt{n}}]}}{ \sqrt{\frac{1}{\sqrt{n}}} } \right] \leq n^{1/4},
\end{align*}
and also, since $\PP(\tfrac{1}{n+1} \leq Q_n \leq 1-\tfrac{1}{n+1})=1$,
\begin{align*}
\EE\left[ \frac{\I{Q_n \notin [\frac{1}{\sqrt{n}}, 1-\frac{1}{\sqrt{n}}]}}{\sqrt{Q_n(1-Q_n)}} \right] \leq \EE\left[ \frac{\I{Q_n \notin [\frac{1}{\sqrt{n}}, 1-\frac{1}{\sqrt{n}}]}}{ \sqrt{1/(n+1)} } \right] = \frac{ \PP \left( Q_n \notin [\tfrac{1}{\sqrt{n}}, 1-\tfrac{1}{\sqrt{n}}] \right) } {1/\sqrt{n+1}}.
\end{align*}
Hence,
\begin{align*}
\frac{1}{\sqrt n} \EE\left[ \frac{1}{\sqrt{Q_n(1-Q_n)}} \right]
& \leq \frac{n^{1/4}}{\sqrt n} + \frac{\sqrt{n+1}}{\sqrt n} \PP \left( Q_n \notin [\tfrac{1}{\sqrt{n}}, 1-\tfrac{1}{\sqrt{n}}] \right) \\
& \leq \frac{1}{n^{1/4}} + 2 \PP \left( Q_n \notin [\tfrac{1}{\sqrt{n}}, 1-\tfrac{1}{\sqrt{n}}] \right).
\end{align*}
Lastly, note that by Hypothesis~\hyperlink{hypD}{(D)}
\begin{align*}
\PP \left( Q_n \notin [\tfrac{1}{\sqrt{n}}, 1-\tfrac{1}{\sqrt{n}}] \right) \leq \int_0^1 g(q) \I{q \notin [\tfrac{1}{\sqrt{n}}, 1-\tfrac{1}{\sqrt{n}}]}~\dd q \to 0
\end{align*}
as $n \to +\infty$.
This proves~\hyperlink{hypB}{(B)}, which concludes the proof.
\end{proof}

We now give the proof of Proposition~\ref{propos}.
It is obtained from Corollary~\ref{corol}.

\begin{proof}[Proof of Proposition~\ref{propos}.]
The proof consists on showing, first, that for any $b>0$ the sequence of signatures $(s^{(n)}_n : n \geq 1)$ defined in Expression~\eqref{propos:sign} satisfies Hypotheses~\hyperlink{hypC}{(C)} and~\hyperlink{hypD}{(D)} for the density $f (q) = b (1-q)^{b-1}$ for $q$ in $(0,1)$; and second, showing that $\int_0^1 \PP(S_t < -\log(1-q)) f_Q(q) \dd q = e^{-\psi(b) t}$ for all $t>0$.

To that end, recall first that the signature of the $n$-th system is $s_k^{(n)} = \left( n-k+1 \right)^{b-1} C^{(n)}$ for $k = 1, \ldots, n$, where $C^{(n)}$ is a normalizing constant such that $\sum_{k=1}^n s_k^{(n)} = 1$.
We will use the following result.
\begin{lemma}\label{lem} 
For all $b>0$ the asymptotic behavior of $C^{(n)}$ is given by
\begin{align*}
\lim_{n\to \infty }C^{(n)} \frac{n^b}{b} = 1.
\end{align*}
\end{lemma}

We prove first Hypothesis~\hyperlink{hypC}{(C)}.
Note that
\begin{align}\label{prop:eq 1}
n ~ s_{\lceil n q \rceil } ^{(n)} = n \left( n+1-\lceil n q \rceil  \right)^{b-1} C^{(n)} = C^{(n)} \frac{n^b}{b} \, b \left( 1 - \frac{\lceil n q \rceil }{n} +\frac{1}{n} \right)^{b-1} .
\end{align}
If $b\neq 1$, clearly $n ~ s_{\lceil n q \rceil } ^{(n)} \to f(q)$ as $n \to +\infty$, by Lemma~\ref{lem} and
\begin{align}\label{eq:bound}
q \leq \frac{\lceil n q \rceil}{n}  \leq q+\frac{1}{n}.
\end{align}
Otherwise, when $b=1$,  we have a discrete uniform distribution with $C^{(n)}=1/n$ and the convergence to the uniform density is straightforward.
This proves~\hyperlink{hypC}{(C)}.
 
We now check the dominated convergence Hypothesis~\hyperlink{hypD}{(D)}.
For that, using Lemma \ref{lem} there exists a constant $K_b$ depending on $b$ but not on $n$ such that $C^{(n)}n^b\leq K_b$ for all $n\geq 1$.
We will consider separately the cases where $b$ in $(0,1)$, $b=1$, and $b>1$.

We first consider the case $b$ in $(0,1)$.
Using in~\eqref{prop:eq 1} that $C^{(n)}n^b\leq K_b$ and $\left( 1 - q + (1-\lceil nq \rceil + nq ) / n \right)^{b-1} \leq (1-q)^{b-1}$ ---from~\eqref{eq:bound} and the fact that $b-1<0$---, we obtain
\begin{align}\label{proof prop ineq1}
n \, s^{(n)}_{\lceil nq \rceil}  \leq K_b (1-q)^{b-1}
\end{align}
for all $n \geq 1$ and $q$ in $(0,1)$, with $\int_0^1 (1-q)^{b-1} \dd q = 1/b < +\infty$.
Next, the case $b=1$ is trivial, since $n \, s^{(n)}_{\lceil nq \rceil}  = 1$ for all $n$ and $q$ in $(0,1)$.
Finally, the case $b>1$ holds from the fact that
\begin{align}\label{proof prop ineq2}
n \, s^{(n)}_{\lceil nq \rceil} \leq K_b b (2-q)^{b-1}
\end{align}
for all $n \geq 1$ and $q$ in $(0,1)$, with $\int_0^1 b (2-q)^{b-1} \dd q =b(2^b-1)< +\infty$.
Indeed, from~\eqref{eq:bound} and since $b>1$ we obtain $\left( 1 - q - (\lceil nq \rceil - nq -1) / n \right)^{b-1} \leq (1-q + 1)^{b-1}$.
This proves~\hyperlink{hypD}{(D)}.

We have thus checked that for all $b >0$ the sequence of signatures $(s^{(n)}_n : n \geq 1)$ defined in Expression~\eqref{propos:sign} satisfies Hypotheses~\hyperlink{hypC}{(C)} and~\hyperlink{hypD}{(D)} for the density $f(q) = b (1-q)^{b-1}$ for $q$ in $(0,1)$.
Hence, by Corollary~\ref{corol}, we obtain that $\PP ( T_\text{sys}^{(n)} > t ) \to \PP(S_t < -\log(1-Q))$ for all $t>0$ as $n \to +\infty$, where $Q$ is a random variable over $(0,1)$ with density $b (1-q)^{b-1}$ for $q$ in $(0,1)$; that is, $Q$ is a beta(1, $b$) random variable. Lastly, note that
\begin{align*}
\PP(S_t < -\log(1-Q)) = \PP(1-Q < e^{-S_t}) = \EE \left[ e^{-b S_t} \right] = e^{-t \psi(b)},
\end{align*}
where the second equality comes from basic properties of the beta distribution, and the last equality comes from the definition of the Laplace exponent $\psi$ of the L\'evy subordinator $S$.

It remains to prove Lemma~\ref{lem}.
The case $b = 1$ is trivial since $C^{(n)}=n^{-1}$ and $C^{(n)}n b=1$.
So let us consider $b\neq 1$.
Note that $x \in \RR_+ \mapsto x^{b-1}$ is a monotonic function in $x$.
For $b>1$ it is nondecreasing, so $\int_1^n (x-1)^{b-1} \dd x \leq \sum_{i=1}^n i^{b-1} \leq \int_1^n x^{b-1} \dd x$; and since $C^{(n)} = 1/\sum_{i=1}^n i^{b-1}$ we get that
\begin{align*}
\frac{b}{ (n-1)^{b}} \geq C^{(n)} \geq \frac{b}{n^b-1},
\end{align*}
so re-arranging terms,
\begin{align*}
\frac{1}{1-1/n^b} \leq C^{(n)} \frac{n^b}{b} \leq \left(\frac{1}{1-1/n}\right)^{b}.
\end{align*}
For $b\in (0,1)$ we use that the function $x \mapsto x^{b-1}$ is nonincreasing, obtaining the analogous bound,
\begin{align*}
 \left(\frac{1}{1-1/n}\right)^{b} \leq C^{(n)}  \frac{n^b}{b}  \leq  \frac{1}{1-1/n^b} .
 \end{align*}
This proves that for all $b> 0$ we have $\lim_{n \to \infty} C^{(n)}  \frac{n^b}{b}=1$, proving Lemma~\ref{lem}.
 This concludes the proof of Proposition~\ref{propos}.
\end{proof}

\section{Conclusions}\label{sec:conclusions}

In this work we have proposed the first ---to the best of the author's knowledge--- tractable asymptotic result on the reliability function of a system whose number of components is growing to infinity.
We do this for systems with a mixed coherent structure, see Section~\ref{subsec:mixed coherent}, and whose components fail at random times following a LFMO  distribution, see Section~\ref{subsec:LFMO}.
Indeed, in the main result of this work, Theorem~\ref{theo} in Section~\ref{sec:main}, we show that in the limit the reliability function, i.e., the probability that the system is still working by a certain given time, converges to the probability of a first-passage time of the underlying L\'evy-subordinator of the LFMO distribution.
In more probabilistic terms, we show that the system failure time converges in distribution to a first-passage time of the subordinator.
For this, we make Hypotheses~\hyperlink{hypA}{(A)} and~\hyperlink{hypB}{(B)} on the signature of the system; in rough terms, they impose that the natural scale at which systemic failure occurs is when a fraction of the components have failed.
In particular, this includes the case where the system fails when a certain given fraction of the total number of components has failed.
To illustrate our approach, in Proposition~\ref{propos} in Section~\ref{sec:example} we give an example where the system failure time converges in distribution to an exponential random variable.

We remark that we have considered a setting where we can separate the failure times of the components from the fact of the system being working or not; indeed, the notion of mixed coherence gives the probability of the system being working or not under a given static configuration of working and non-working components, and such probability does not depend on the time when the non-working components failed, nor on how long have the working components been working.
In other words, we have separated the analysis of the system into a \emph{static} part where there is no time involved and where the status of the system only depends on the status of its individual components, and a \emph{dynamic} part where we only focus on the status of the components throughout time.
Our computational results show that our limit result can approximate well the reliability function even when the system has only a few hundred components.

\subsection{Open questions}

An interesting question that our work opens is how can one qualify that a system satisfies Hypotheses~\hyperlink{hypA}{(A)} and~\hyperlink{hypB}{(B)}, or even a given signature such as the one we give in Proposition~\ref{propos}?
Especially because computing exactly the signature of a mixed system can be challenging, see, e.g.,~Definition~\ref{def:signature}.
In the remarks following Proposition~\ref{propos} we have mentioned that the qualification can be done heuristically, based on statistical notions on which failures tend to take down the system.
Another possibility is using a probabilistic decision scheme where the system in use is randomized and chosen as the system $k$-out-of-$n$ with probabilities given by the desired system signature ---see the remarks following Definition~\ref{def:signature}---, however this can be non-practical, if not artificial.
Nonetheless, it can be of interest to derive a more formal or general method to check that the system under study satisfied the hypotheses made on the system signature.

Another intriguing question that our work opens is if the mean time-to-failure of the system, i.e., the mean of the system failure time, converges to the mean of the limit first-passage time.
In technical terms, this is a question on interchange of limit and expectation, and it would establish a mode of convergence that is stronger than the convergence in distribution we have shown here.
This is an important question since the mean time-to-failure is of high interest in Reliability Engineering, and our limit can open the door to tackling it by simulation and analysis of the first-passage time of the underlying L\'evy-subordinator.

Finally, in the setting where components fail according to the LFMO distribution, the following converse problem can be of interest: what are the sufficient conditions on the system signature for asserting strong unimodality of the system lifetime? This exciting question has recently been addressed in~\cite{rychlik2021signature}  when components have iid exponential distributed lifetimes. Our results could allow to answer this question when components have exchangeable lifetimes distributions.



\bmsection*{Acknowledgments}
JB acknowledges the financial support of FONDECYT grants 1161064 and 1231207, and of Programa Iniciativa Cient\'ifica Milenio NC120062.
Part of this paper was written when JB was Visiting Associate Professor \&  Esbach Scholar at the Industrial Engineering and Management Sciences Department at Northwestern University.
GL acknowledges the financial support of FONDECYT grant 11230256 and of CORFO InnovaChile project 14ENI2-26865.
PR acknowledges the financial support of project FCE-ANII grant 156693.


\bmsection*{Conflict of interest}
The authors declare no potential conflict of interests.

\bibliography{wileyNJD-AMS}


\end{document}